\documentstyle[10pt,amsfonts,german]{article}

\newtheorem{theorem}{Theorem}[section]
\newtheorem{satz}[theorem]{Theorem}  
\newtheorem{lem}[theorem]{Lemma}  
\newtheorem{cor}[theorem]{Corollary}  
\newtheorem{defi}[theorem]{Definition}  
\newcommand{\titel}[1]{\begin{center} {\LARGE\bf\sc\ #1}
\end{center}\title{#1}\vspace{1cm}}
\newcommand{\autor}[1]{\begin{center}{\sc #1}\end{center} \author{#1}}
\newcommand{\adresse}[1]{ \begin{center} {\it #1}
\end{center}\vspace{0,3cm}}
\newcommand{\abstrakt}[4]{\mbox{} \hrulefill \mbox{} 
\begin{quotation} 
\noindent
{\footnotesize\bf Abstract: }{\footnotesize #1 }\\[0,5cm] 
 {\footnotesize\it Subj. Class:} {\footnotesize #2}\\
  {\footnotesize\it 1991 MSC:} {\footnotesize #3} \\ {\footnotesize\it Keywords:}
 {\footnotesize #4}
\end{quotation} 
 \vspace{0,1cm}
\mbox{} \hrulefill \mbox{}\vspace{1cm}}
\newcommand{\map}[3]{#1:#2\longrightarrow #3}  

\newcommand{\krumm}{{\cal R}}

\newcommand{\r}{ {\mathbb{R} } }
\newcommand{\com}{ { \mathbb{C}} }

\newcommand{\beweis}{ {\sc Proof:} }
\newcommand{\bemerkung}{ {\bf Remark:} }

\newcommand{\schluss}{\nopagebreak  \hfill $\Box$ \\[0,2cm] }  

\newcommand{\glungen}[1]{\begin{eqnarray*}#1 \end{eqnarray*}}
\newcommand{\lab}[1]{\label{#1}}
\newcommand{\glei}[2]{\begin{equation}\lab{#1} #2 \end{equation}}
\newcommand{\dichte}{ {\cal L}}

\newcommand{\enumer}[1]{
                        \renewcommand{\labelenumi}{\arabic{enumi}.}
                        \begin{enumerate}#1\end{enumerate}
                        }

\catcode`\"=\active
\catcode`\=\active
\catcode`\"=\active
\catcode`\Ž=\active
\catcode`\š=\active
\catcode`\™=\active
\catcode`\á=\active
\def"{\"a}
\def{\"u}
\def"{\"o}
\defŽ{\"A}
\defš{\"U}
\def™{\"O}
\letá=\ss
\newcommand{\sym}{ \textsl{sym} }

\newcommand{\Alt}{ \textsl{ Alt} }

\newcommand{\einhalb}{\frac{1}{2}}  

\begin{document}

\titel{a note on real killing spinors in weyl geometry \footnote{Supported by the SFB 288 of the DFG}}
\autor{Volker Buchholz\footnote{ e-mail: {\tt bv@mathematik.hu-berlin.de}}}
\adresse{ Humboldt Universit"at zu Berlin, Institut f"ur Reine Mathematik,\\ Ziegelstra"se 13a , D-10099 Berlin. } 
\abstrakt{This text is dedicated to the real Killing equation on 3-dimensional
Weyl manifolds. Any manifold admitting a real Killing spinor of weight 0
satisfies the conditions of a Gauduchon-Tod geometry. Conversely,
any simply connected Gauduchon-Tod geometry has a 2-dimensional space of
solutions
of the real Killing equation on the spinor bundle of weight 0.}{Differential Geometry}{53C05;53C10;53A30}
{ Weyl geometry; Spin geometry; Killing equation}

\section{Introduction}

In \cite{sweyl} we introduced the spinor geometry on Weyl
manifolds and investigated the Dirac-, Twistor- and Killing equation in this
context. Concerning the real Killing equation we presented in \cite{sweyl} the following result:
\begin{satz}[see \cite{sweyl}, Theorem 3.1]
   Let $\psi\in\Gamma(S^w)$ be a real Killing spinor on a Weyl manifold
	$(M^n,c,W)$, i.e. there exists a
	\lab{killingsatz}
	density $\beta \in \Gamma(\dichte^{-1})$ for which 
	$$
   \nabla^{S,w}\psi=\beta\otimes\nu\psi, \quad \beta\in \Gamma(\com \otimes
   \dichte^{-1})
   $$
 	is satisfied. Then the following statements hold:  
	\enumer
		{
        \item $R = 4n(n-1)\beta^2$. 
                \item $w \neq 0$: $W$ is exact and Einstein-Weyl.
                \item $w=0$, $n \geq 4$: $W$ is exact and Einstein-Weyl.
		}
\end{satz}
The following equations were obtained within the proof and will play an important
role in the sequel:
\glei{grundgleichung}
		{
        \mu^2Ric'\otimes \psi  =  2\left(n-1-\frac{n-1}{n-2}\right)\nabla\beta\otimes\psi+\left(1-\frac{n-1}{n-2}\right)
        \mu^{12}\Alt\nabla\beta\otimes c\otimes\psi +\frac{R}{n}\nu\psi-\mu^2
		  F\otimes\psi
		}
\glei{nunablabeta}
	{
        F\cdot\psi = -\frac{4(n-1)}{n-2}\nabla\beta\cdot\psi.
	}		  
Theorem \ref{killingsatz} gives no statement for the case $n=3$ and $w=0$. In the next
section we proof that in three dimensions the existence of a real Killing spinor of weight $0$ is
essentially equivalent to the fact that this manifold is a Gauduchon-Tod
geometry:
\begin{defi}[see \cite{gautod}, Proposition 5] \label{defgautod}
A $3$-dimensional Weyl manifold $(M^3,W,c)$ is called Gauduchon-Tod geometry, if  
there exists a density $\beta \in \Gamma(\dichte^{-1})$ such that the following
conditions are satisfied:
\enumer
        {
		  \item $W$ is Einstein-Weyl;
        \item $R=24\beta^2$;
        \item $4 \nabla\beta =*F$.
        }
\end{defi}

\bemerkung The $\kappa\in C^\infty(M,\r)$ in (\cite{gautod},  Proposition 5) is related to
the $\beta\in \Gamma(\dichte^{-1})$ of Definition \ref{defgautod} in the following way: 
$
\kappa l^{-1}_{g_\Sigma}=-4\beta. 
$
For more
information on Gauduchon-Tod geometries, e.g. their classification,  see \cite{gautod}  and the  references
therein.

\bigskip

Hence, the main result of this text is as follows:
\begin{satz} \label{hauptkillingsatz}
Let $(M^3,c,W)$ be a $CSpin$-manifold.
\enumer
	{ 
	\item If $\psi\in\Gamma(S^0)$ is a real Killing spinor then the
	space of solutions of the Killing equation is 2-dimensional and $(M^3,c,W)$ is 
	a Gauduchon-Tod geometry. 
	\item Conversely, any simply connected Gauduchon-Tod geometry has a two
	dimensional space of Killing spinors of weight $0$.      
	}	
\end{satz}

\bigskip

{\bf Acknowledgement:}
I would like to thank N. Ginoux for pointing out a gap in a former version of
\cite{sweyl}. This concerns the real Killing equation in dimension $3$ and for
weight $0$.

\section{The proof of Theorem \ref{hauptkillingsatz}}

Let $(M^3,c,W)$ be a CSpin-manifold. The curvature tensor of the Weyl structure $W$ is given by 
$$
\krumm= Ric^N \bigtriangleup c + F\otimes c,
$$
where 
$\map{\bigtriangleup}{T^{2,0} \times T^{2,0}}{T^{4,0}}$ 
$$
	\omega\bigtriangleup
	\eta:=[(23)+(12)(24)(34)-(24)-(12)(23)]\omega\otimes\eta, \quad \omega,\eta
	\in T^{2,0}
$$
is the so called Kulkarni-Nomizu product (see \cite{besse}) and 
\glei{normricci}
	{
	Ric^N:=-{\sym}_0Ric-\frac{1}{12}Rc+\frac{1}{2}F
	}
is the normalized Ricci tensor of $W$ (see \cite{peca}). $\sym_0$ denotes the
symmetric trace free  part of a $(2,0)$-tensor. The following Lemma is a tool
for calculations with Kulkarni-Nomizu products in spin geometry:
\begin{lem} \label{dreieckmult} Let $\omega$ be a $(2,0)$-tensor. Then the
following algebraic identity holds in any dimension:
$$
\mu^{34}\omega\bigtriangleup c=2 \Alt\nu\mu^2\omega-2\Alt\omega.
$$	
\end{lem}
\beweis
\glungen
	{
	\mu^{34}\omega\bigtriangleup c &=& \mu^{34}[(23)+(12)(24)(34)-(24)-(12)(23)]
	\omega\otimes c
	= [\mu^{24}+(12)\mu^{32}-\mu^{32}-(12)\mu^{24}]\omega\otimes c\\
	&=& [\mu^{23}+(12)\mu^{32}-\mu^{32}-(12)\mu^{23}]\omega\otimes c
	= [-2\mu^{32}+2tr^{23}+2(12)\mu^{32}-2(12)tr^{23}]\omega\otimes c\\
	&=& -2[\mu^{32}-(12)\mu^{32}]\omega\otimes c - 2
	{\Alt}\omega
	= -2[1-(12)]\nu\mu^2\omega\otimes c - 2
	{\Alt}\omega\\
	&=&2 \Alt\nu\mu^2\omega-2\Alt\omega.	
	}
\schluss

\begin{lem}
Let  $4 \nabla\beta = *F$ be satisfied on $(M^3,c,W)$. Then the following
identities are true for any spinor $\psi \in \Gamma(S^w)$:
\glei{ersteslemma}
	{
	\frac{1}{4}\Alt\nu\mu^2F\otimes\psi - \Alt\nabla\beta\otimes\nu \psi - \einhalb
	F\otimes \psi  =  0.
	}
and 
\glei{zweiteslemma}
	{	
	(\nu\nabla\beta-\nabla\beta\cdot\nu)\cdot\psi -\einhalb\mu^2 F\otimes\psi 
	=  0.
	}
\end{lem}

\beweis Denote  by $(e_1,e_2,e_3)$ a local weightless conformal frame 
 on $(M^3,c)$ as well as $(\sigma_1,\sigma_2,\sigma_3)$ its dual. 
In dimension 3 we have the important relation
\glei{dim3rech}
	{
	 e_i\cdot e_j \cdot \psi = -\sum^3_{k=1}\epsilon_{ijk}e_k\cdot\psi.
	}
Here $\epsilon_{ijk}$ denotes the Levi-Civita symbol. 
Since the $*$-operator on 2-forms is defined by the formula
$$
*F=\einhalb\sum^3_{i,j=1}F(e_i,e_j)*(\sigma_i\wedge\sigma_j)=\einhalb\sum^3_{i,j,k=1}F(e_i,e_j)\epsilon_{ijk}\sigma_k
$$
we can rewrite the assumption as follows: 
\glei{monop}
	{
	8\nabla_{e_k}\beta=\sum^3_{i,j=1}\epsilon_{ijk}F(e_i,e_j).
 	}
From the algebraic identity
$$
	F\cdot\psi  =  \sum^3_{i,j=1}F(e_i,e_j)e_i\cdot e_j\cdot \psi
	 =  -\sum^3_{i,j,k=1}\epsilon_{ijk}F(e_i,e_j) e_k\cdot \psi
	 =  -8\sum^3_{k=1}(\nabla_{e_k}\beta) e_k\cdot \psi
$$	
and (\ref{monop}) then follows:
\glei{monopnichthodge3}
	{
	F\cdot\psi = -8\nabla\beta\cdot\psi,
	}	
Using
(\ref{monop}) and (\ref{dim3rech})  we get:

\begin{eqnarray*}
\lefteqn{\frac{1}{4}Alt\nu\mu^2F\otimes\psi - Alt\nabla\beta\otimes\nu \psi - \einhalb
F\otimes \psi} \\
& = &
\frac{1}{4}\sum^3_{i,j,k=1} F(e_j,e_k) \sigma_i\wedge\sigma_j \otimes e_i \cdot
e_k\cdot\psi-
\sum^3_{i,j=1}(\nabla_{e_i}\beta) \sigma_i\wedge\sigma_j \otimes e_j\cdot\psi
- \frac{1}{4}\sum^3_{i,j=1} F(e_i,e_j)\sigma_i\wedge\sigma_j\otimes \psi\\
& = &-\frac{1}{4}\sum^3_{i,j=1} F(e_j,e_i) \sigma_i\wedge\sigma_j \otimes \psi
+\frac{1}{4}\sum^3_{\stackrel{i,j,k=1}{i\neq k}} F(e_j,e_k)
\sigma_i\wedge\sigma_j \otimes e_i \cdot e_k\cdot\psi\\
&  & -\sum^3_{i,j=1}(\nabla_{e_i}\beta)\sigma_i\wedge\sigma_j \otimes e_j\cdot\psi -
\frac{1}{4}\sum^3_{i,j=1} F(e_i,e_j)\sigma_i\wedge\sigma_j\otimes \psi\\
& = &   \frac{1}{4}\sum^3_{i,j,k,l=1}\epsilon_{kil} F(e_j,e_k)
\sigma_i\wedge\sigma_j \otimes
e_l\cdot\psi-\frac{1}{8}\sum^3_{i,j,k,l=1}\epsilon_{kli}F(e_k,e_l)
\sigma_i\wedge\sigma_j \otimes e_j\cdot\psi\\
& = &   \frac{1}{4}\sum^3_{i,j,k=1}\epsilon_{kij} F(e_j,e_k)
\sigma_i\wedge\sigma_j \otimes
e_j\cdot\psi-\frac{1}{8}\sum^3_{\stackrel{i,j,k,l=1}{j=k}}\epsilon_{jli}F(e_j,e_l)
\sigma_i\wedge\sigma_j \otimes
e_j\cdot\psi\\
& - &\frac{1}{8}\sum^3_{\stackrel{i,j,k,l=1}{j=l}}
\epsilon_{kji}F(e_k,e_j)\sigma_i\wedge\sigma_j \otimes e_j\cdot\psi\\
& = &   \frac{1}{4}\sum^3_{i,j,k=1}\epsilon_{kij} F(e_j,e_k)
\sigma_i\wedge\sigma_j \otimes
e_j\cdot\psi-\frac{1}{8}\sum^3_{i,j,l=1}\epsilon_{lij}F(e_j,e_l)
\sigma_i\wedge\sigma_j \otimes e_j\cdot\psi\\
& - & \frac{1}{8}\sum^3_{i,j,k=1}
\epsilon_{kij}F(e_j,e_k)\sigma_i\wedge\sigma_j \otimes e_j\cdot\psi\\
& = &0,
\end{eqnarray*}
i.e. we have shown (\ref{ersteslemma}). 
Contracting  (\ref{ersteslemma}) by  
Clifford multiplication, using (\ref{monopnichthodge3}) and
$\mu^2\nu=-\nu\mu^1-2Id$ yields  :
\glungen
	{
	0 & = &	\frac{1}{4}\mu^2Alt\nu\mu^2F\otimes\psi - 
	\mu^2Alt\nabla\beta\otimes\nu \psi - \einhalb\mu^2F\otimes \psi    \\
	& = & \frac{1}{4}\mu^2\nu\mu^2F\otimes\psi - 
	\mu^2\nabla\beta\otimes\nu \psi - \einhalb\mu^2F\otimes \psi
	-\frac{1}{4}\mu^1\nu\mu^2F\otimes\psi + 
	\mu^1\nabla\beta\otimes\nu \psi \\
	& = & -\frac{1}{4}\nu F\cdot\psi -\einhalb \mu^2F\otimes \psi  
	+3\nabla\beta\otimes \psi - \einhalb\mu^2F\otimes \psi
	+\frac{3}{4}\mu^2F\otimes\psi + 
	\nabla\beta\cdot\nu \psi\\ 
	& = & -\frac{1}{4}\nu F\cdot\psi -\einhalb \mu^2F\otimes \psi  
	+3\nabla\beta\otimes \psi - \einhalb\mu^2F\otimes \psi
	+\frac{3}{4}\mu^2F\otimes\psi + 
	\nabla\beta\cdot\nu \psi \\
	& = & \nu \nabla\beta\cdot\psi   
	+\nabla\beta\otimes \psi -\frac{1}{4}\mu^2F\otimes\psi\\
	& = & \einhalb (\nu
	\nabla\beta-\nabla\beta\cdot\nu)\cdot\psi-\frac{1}{4}\mu^2F\otimes\psi. 
	}
Hence (\ref{zweiteslemma}) is true.
\schluss
After these preliminary calculations
we get to the 
 proof of Theorem \ref{hauptkillingsatz}:
 Let $\psi\in\Gamma(S^0)$ be a real Killing spinor. The first
	statement follows immediately from the existence of an equivariant quaternionic 
	structure $\bf{j}$ on the spinor module, which commutes with the Clifford
	multiplication. By Theorem 
	\ref{killingsatz} and its proof we have 
	$$
	 	R=24\beta^2;\quad F\cdot\psi = -8\nabla\beta\cdot\psi.
	$$
	Since $\psi$ vanishes nowhere the second equation is equivalent to
	$$
	 	         4 \nabla\beta =*F
	$$				
	by (\ref{dim3rech}).
	Therefore we have already verified the conditions {\em 2.} and {\em 3.} of
	the definition of a Gauduchon-Tod geometry. It remains to proof that the
	manifold is Einstein-Weyl. To this end we have to simplify 
	(\ref{grundgleichung}) by means of $R=24\beta^2$ and $Ric'=\sym_0 Ric
	+\frac{1}{3}Rc-\einhalb F$ to
	$$
	\mu^2\sym_0Ric'\otimes\psi =
							-(\nabla\beta\cdot\nu-\nu\nabla\beta)\cdot\psi -
							\einhalb\mu^2 F\otimes\psi	
	$$
	But the righthandside vanishes according to (\ref{zweiteslemma}). Hence $W$
	is Einstein-Weyl. 
	\bigskip
	
	Conversely, let $W$ be a simply connected Gauduchon-Tod geometry. Is is
	sufficient to show that $S^0$ is flat with
	respect to $\nabla^\beta=\nabla^{S,0} - \beta \otimes\nu$. To this end, we have to
	proof that the curvature of $\nabla^\beta$ vanishes. We use the properties of
	Gauduchon-Tod geometries given in Definition \ref{defgautod}, the result of
	 Lemma
	\ref{dreieckmult} and the equations  (\ref{ersteslemma}) 
	and $\krumm^{S,0}=\frac{1}{4}\mu^{34}Ric^N\bigtriangleup c=-
	\frac{1}{4}\mu^{34}\left(\frac{1}{12}Rc-\einhalb F\right)\bigtriangleup
	c$.
\glungen
	{
	\krumm^\beta &=& \Alt \nabla^\beta\circ\nabla^\beta
	 = 
	\Alt
	\nabla^\beta\circ (\nabla^{S,0} - \beta\otimes \nu)\\
	&=& \Alt\left( \nabla^{S,0}\circ\nabla^{S,0}
	-\nabla\beta\otimes\nu-(12)\beta\otimes\nu\nabla^{S,0}-
	\beta\nu\nabla^{S,0}+\beta^2\nu\nu\right)\\
	&=& \krumm^{S,0}- \Alt (\nabla\beta)\nu
	+\beta^2\Alt\nu^{12}
	= - \frac{1}{4}\mu^{34}\left(\frac{1}{12}Rc-\einhalb F\right)\bigtriangleup c -\Alt (\nabla\beta)\nu
	+\beta^2\Alt\nu^{12}\\
	&=&	 -
	\frac{1}{4}\left(\frac{1}{6}R\Alt\nu \mu^2c- \Alt\nu\mu^2 F+ 2 F \right) -\Alt (\nabla\beta)\nu
	+\beta^2\Alt\nu^{12} \\
	& = & -\frac{1}{24}R\Alt\nu^{12} + 	\frac{1}{4}\Alt\nu\mu^2 F-\einhalb F 
	-\Alt (\nabla\beta)\nu +\beta^2\Alt\nu^{12} \\
	& = & \frac{1}{4}\Alt\nu\mu^2 F-\einhalb F 
	-\Alt (\nabla\beta)\nu\\
	& = & 0.
	} 
\schluss		
\selectlanguage{\english}


\begin{thebibliography}{1}

\bibitem{besse}
A.L. Besse.
\newblock {\em Einstein manifolds}.
\newblock Springer Verlag, 1987.

\bibitem{sweyl}
V.~Buchholz.
\newblock Spinor equations in {W}eyl geometry.
\newblock {\em {\em to appear in} Proc. of the Winter School on Geometry and
  Physics, Srni}, 1999; math.DG/9901125 v2.

\bibitem{peca}
D.~Calderbank{,}~H. Pederson.
\newblock Einstein-{W}eyl geometry.
\newblock {\em Odense Universitet, preprint}, No. 40, 1997.

\bibitem{gautod}
P.~Gauduchon{,}~K.P. Tod.
\newblock Hyper-hermitian metrics with symmetry.
\newblock {\em Journal of Geometry and Physics}, pages 291--304, 1998.

\end{thebibliography}
\end{document}